\documentclass[12pt,reqno]{amsart}
\title{Levi-Civita connections and vector fields \\ ~ \\ for noncommutative differential calculi}
\date{v2 April 2020}
\usepackage{latexsym, amsfonts, amsmath, amssymb}
\usepackage{fullpage}

\usepackage{color}

\newcommand{\beq}{\begin{equation}}
\newcommand{\eeq}{\end{equation}}

\newcommand{\be}{\begin{equation}}
\newcommand{\ee}{\end{equation}}
\newcommand{\bea}{\begin{eqnarray}}
\newcommand{\eea}{\end{eqnarray}}
\newcommand{\bean}{\begin{eqnarray*}}
\newcommand{\eean}{\end{eqnarray*}}

\newcommand{\brray}{\begin{array}}
\newcommand{\erray}{\end{array}}
\newcommand{\ben}{\begin{equation}{nonumber}}
\newcommand{\een}{\end{equation}{nonumber}}

\newtheorem{defn}{Definition}[section]
\newtheorem{thm}[defn]{Theorem}
\newtheorem{lemma}[defn]{Lemma}
\newtheorem{prop}[defn]{Proposition}
\newtheorem{corr}[defn]{Corollary}
\newtheorem{xmpl}[defn]{Example}
\newtheorem{rmk}[defn]{Remark}
\newcommand{\bdfn}{\begin{defn}}
\newcommand{\bthm}{\begin{thm}}
\newcommand{\blmma}{\begin{lemma}}
\newcommand{\bppsn}{\begin{prop}}
\newcommand{\bcrlre}{\begin{corr}}
\newcommand{\bxmpl}{\begin{xmpl}}
\newcommand{\brmrk}{\begin{rmk}}
\newcommand{\edfn}{\end{defn}}
\newcommand{\ethm}{\end{thm}}
\newcommand{\elmma}{\end{lemma}}
\newcommand{\eppsn}{\end{prop}}
\newcommand{\ecrlre}{\end{corr}}
\newcommand{\exmpl}{\end{xmpl}}
\newcommand{\ermrk}{\end{rmk}}
\newcommand{\dubois}{\Omega^1_{{\rm Der}} ( \A )}

\newcommand{\IC}{\mathbb{C}}

\newcommand{\id}{\mathrm{id}}

\newcommand{\oneform}{{{\Omega}^1 ( \mathcal{A} ) }}

\newcommand{\twoform}{{{\Omega}^2}( \mathcal{A} )}
\newcommand{\tensora}{\otimes_{\mathcal{A}}}

\newcommand{\tensorc}{\otimes_{\mathbb{C}}}
\newcommand{\ot}{\otimes}
\newcommand{\A}{\mathcal{A}}

\newcommand{\E}{\mathcal{E}}

\newcommand{\F}{\mathcal{F}}

\newcommand{\Acenter}{\mathcal{Z}( \mathcal{A} )}

\newcommand{\Ecenter}{\mathcal{Z}( \mathcal{E} )}

\newcommand{\Psym}{P_{\rm sym}}
\newcommand{\Hom}{{\rm Hom}}

\newcommand{\omegaone}{\omega_{(1)}}

\newcommand{\omegazero}{\omega_{(0)}}
\newcommand{\etazero}{\eta_{(0)}}
\newcommand{\etaone}{\eta_{(1)}}
\newcommand{\zeroomega}{{}_{(0)}\omega}
\newcommand{\oneomega}{{}_{(1)}\omega}

\newcommand{\vgtwo}{V_{g^{(2)}}}

\newcommand{\RNum}[1]{\uppercase\expandafter{\romannumeral #1\relax}}

\begin{document}
\maketitle
\begin{center}
{\large {Jyotishman Bhowmick, Debashish Goswami, Giovanni Landi}}\\
\end{center}
\begin{abstract}
We study covariant derivatives on a class of centered bimodules 
$\mathcal{E}$ 
over an algebra $\A.$ We begin by identifying a $\Acenter$-submodule $ \mathcal{X} ( \A ) $ 
which can be viewed as the analogue of vector fields in this context; $ \mathcal{X} ( \A ) $ is proven to be a Lie algebra. Connections on $\mathcal{E}$ are in one to one correspondence with covariant derivatives on $ \mathcal{X} ( \A )$. 
We recover the classical formulas of torsion and metric compatibility of a connection in the covariant derivative form. As a result, a Koszul formula for the Levi-Civita connection is also derived. 
\end{abstract}

\vspace{4mm}



\vspace{4mm}

\vspace{0.5in}
\tableofcontents
\parskip = .75 ex

\section{Introduction}	

The notion of Levi-Civita connections and associated curvature formulas  in noncommutative geometry have attracted a lot of attention in recent years. Connections can be viewed both on the level of vector fields or that of forms. Consequently, formulations and existence-uniqueness questions of Levi-Civita connections in noncommutative geometry were made at the level of forms as well as derivations. These include the papers by Dubois-Violette and Michor \cite{dubois}, \cite{dubois2}, the papers \cite{LaMa88}, \cite{{LNW94}}, \cite{DHLS96} and the ones by Rosenberg \cite{rosenberg}, Peterka and Sheu \cite{sheu} and more recently the results by Arnlind et al in \cite{pseudo}, \cite{cylinder}, \cite{tiger} (and references therein). Alternative approaches were taken by Fr\"ohlich et al \cite{frolich} and Heckenberger-Schmuedgen \cite{heckenberger},  as well as by Beggs and Majid and collaborators (see \cite{beggsmajidbook} for a comprehensive account). Finally, yet another approach to Levi-Civita connections (among other things) of working in the set up of braided derivations has been pursued by Weber in  the paper \cite{weber}. 

The present paper is a companion article to our work in \cite{article3} where we identified a set of sufficient conditions (also see \cite{article1}) on a differential calculus $(\Omega^\bullet (\A), d )$ on an algebra $\A$, so that there exists a unique torsionless connection which is compatible with a given pseudo-Riemannian bilinear metric $g$ on the space of one-forms $ \oneform $. 
Following \cite{article1}, we will use the terminology ``tame differential calculus" for a differential calculus satisfying these sufficient conditions (see Definition \ref{tame}). 
In \cite{article3}, we have worked with right-connections on forms coming from the differential calculus, arriving 
at a Koszul-type formula on the level of forms. 
We also showed that our Levi-Civita connection is a bimodule connection in the sense of \cite{beggsmajidbook}. 
We refer to the relevant portions of \cite{article1}, \cite{article3} and \cite{soumalya}  for examples of tame differential calculi.

It is a natural question to ask whether the analysis made in \cite{article3} can also be done at the dual level of derivations. We answer this question in the affirmative in this article. Given a tame differential calculus on $\A$,
and a pseudo-Riemannian bilinear metric $g$ on the space of one-forms $ \oneform $, 
we show that the set 
$$ \mathcal{X} ( \A ):= \{ g ( \omega \tensora - ): \omega \in \oneform \} \subseteq \Hom_\A ( \oneform, \A )  $$
 plays the role of vector fields just as in classical differential geometry. Indeed, $\mathcal{X} ( \A ) $ is a bimodule over the algebra $ \Acenter $ and elements of $ \mathcal{X} ( \A ) $ act as derivations on $\A$ in the following way: 
if $ X \in \mathcal{X} ( \A ) $, then  the map $ \delta_X $ defined by the formula
	$$ \delta_X : \A \rightarrow \A, \quad \delta_X ( a ) =  X ( da ) \, 
	$$
is a derivation on $\A$.    Conversely, if $ \phi $ is an element of $  \Hom_\A ( \oneform, \A ) $, then 
$ \delta_\phi  $  is a derivation of $ \A $ if and only if $ \phi \in \mathcal{X} ( \A ).$  Moreover, $  \mathcal{X} ( \A ) $   is naturally a Lie subalgebra of the set of all derivations from $\A$ to $\A.$  The $\Acenter$-$\Acenter$-bimodule $ \mathcal{X} ( \A ) $   will allow us to define the covariant derivative from a connection and then to write a Koszul formula for the Levi-Civita connection in the covariant derivative formulation. As a byproduct, we will recover the classical equations of torsion and metric-compatibility  of a connection in this set-up. 

Let us describe the plan of the article. We will always be working under the hypothesis that our differential calculus 
$(\E= \oneform,  d) $ is tame. In  Section~\ref{preliminaries}, we set up our notations and then recall the main result of \cite{article3}. In  Section~\ref{section3}, we construct a $\Acenter$-bimodule of derivations $ \mathcal{X} ( \A ) $ on $\A$ and prove that $ \mathcal{X} ( \A ) $ is a Lie algebra. In  Section~\ref{section4}, we show that connections on $\E$ and covariant derivatives on $ \mathcal{X} ( \A ) $ are in one to one correspondence. Finally, in Section~\ref{section5}, we derive the Koszul formula for the covariant derivative of the Levi-Civita connection.

\section{Preliminaries} \label{preliminaries}

We begin by spelling out notations and basic results that we use later on. 
An unadorned tensor product $\ot$ will stand for the tensor product $\tensorc$ over the field of complex numbers. Throughout the article $\mathcal{A} $ will denote a complex algebra and $ \Acenter $ will denote its center. 
A subset $ S  $ of a right $ \A $-module $ \E $ will be called  right $ \A $-total in $ \E $ if   the right $\A$-linear span of $ S  $ equals $ \mathcal{E}$. For $\A$-$\A$-bimodules $ \E $ and $ \F $, the symbol $ \Hom_\A  ( \E, \F ) $ will denote the set of all right  $ \A $-linear maps from $ \E $ to $ \F $. The symbol $\E^*$ will stand for $ \Hom_\A ( \E, \A )$.
For $\A$-$\A$-bimodules $ \mathcal{F} $ and $ \mathcal{F^\prime} $, the bimodule multiplications on $ \Hom_\A ( \F, \F^\prime ) $ is given by: 
\begin{equation} \label{7thdec20196} ( a T ) ( f ) = a T ( f ) \in \F^\prime, \quad (T a) ( f ) = T ( a f ),  
\quad a \in \A, f \in \F, \,\, T \in \Hom_\A ( \F, \F^\prime ). \end{equation}

A differential calculus over $\A$ is the datum $ ( \Omega^\bullet ( \A ), \wedge, d ) $ where $ \Omega^\bullet ( \A ) $ is a direct sum of $\A$-$\A$-bimodules $ \Omega^j ( \A ) $, with $ \Omega^0 ( \A ) = \A$.    The map $ \wedge :  \Omega^\bullet ( \A ) \tensora \Omega^\bullet ( \A ) \rightarrow \Omega ( \A ) $ is an $\A$-$\A$-bimodule map such that 
$ \wedge ( \Omega^j ( \A ) \tensora \Omega^k ( \A )  ) \subseteq \Omega^{j + k} ( \A )$.     Finally, $d$ is a map from $ \Omega^j ( \A ) $ to $ \Omega^{j + 1} ( \A ) $ such that 
				$$ d^2 = 0 \quad \textup{and} \quad ( \omega \wedge \eta ) = d \omega \wedge \eta + ( - 1 )^{{\rm deg} ( \omega ) } \omega \wedge d \eta. $$
				Moreover, we will also assume that  $ \Omega^j ( \A ) $ is the right $ \A $-linear span of elements of the form $ d a_1 \wedge d a_1 \wedge \cdots \wedge d a_{j}. $
	
	Throughout this article $ \E $ will stand for the space of one-forms $\oneform$ of a differential calculus. It will also be assumed that $ \E $ is a finitely generated projective right $  \A  $-module.

		\bdfn \label{rLr}
Let $ ( \Omega^\bullet ( \A ), d ) $ be a differential calculus on $\A$.     A (right) connection on $\E: = \oneform$
is a linear map $\nabla :\E  \rightarrow \E \tensora \oneform$ satisfying the Leibniz rule 
$$ \nabla(\omega a)=\nabla(\omega)a + \omega \tensora da$$
for all $ \omega \in \E, a \in \A. $

	The torsion of a connection $ \nabla: \E \rightarrow \E \tensora \oneform $ is the right $\A$-linear map  
	$$ T_\nabla:= \wedge \circ \nabla + d \, : \, \E \rightarrow \twoform. $$
	A connection $\nabla$  is called torsion-less if $ T_\nabla = 0.$
	\edfn
	
	\subsection{Tame differential calculi}
	In this article we will deal with a special type of differential calculus which, following the terminology in \cite{article1} we call \emph{tame} differential calculus. Let us recall that if $\E$ is an $\A$-$\A$-bimodule, the center of $\E$,  denoted by $\Ecenter$ is the set
	$$ \mathcal{Z} ( \mathcal{E} ) = \{ e \in \mathcal{E}: e a = a e ~ \forall ~ a ~ \in \A  \}. $$
	It is easy to check that $ \Ecenter$ is a $\Acenter$-$\Acenter$-bimodule.
	
	\bdfn \label{tame}
A differential calculus $ ( \Omega^\bullet ( \A ), d ) $ on  $ \A $ is called tame if the following conditions are satisfied:
\begin{itemize}
\item[1.] The space of one-forms is given as $ \oneform:=\E = \Ecenter \otimes_{\Acenter} \A$. 

\item[2.] The following short exact sequence of right $\A$-modules splits:
 $$ 
 0 \rightarrow {\rm Ker} ( \wedge ) \rightarrow \E \tensora \E \rightarrow {\rm Ran} ( \wedge ) = \twoform \rightarrow 0. 
 $$ 

\item[3.] Suppose $\sigma = 2 \Psym - 1$,  where $ \Psym $ is the idempotent in $ \Hom_\A ( \E \tensora \E, \E \tensora \E ) $ with image $ {\rm Ker} ( \wedge ) $ and kernel $ \F $ (the complement of ${\rm Ker} ( \wedge )$ in $\E \tensora \E $). Then $ \sigma ( \omega \tensora \eta ) = \eta \tensora \omega $ for all $ \omega, \eta \in \Ecenter. $
\end{itemize}
\edfn
\noindent
A few remarks are in order. Firstly, 
part 2. of Definition \ref{tame} implies that
$$ \E \tensora \E = {\rm Ker} ( \wedge ) \oplus \F, $$
where $\F$ is a right $\A$-module such that $ \wedge: \F \rightarrow \twoform $ is a right $\A$-linear isomorphism. 
In turn, this implies the existence of an idempotent $\Psym$ in $ \Hom_\A ( \E \tensora \E, \E \tensora \E ) $ with  $ {\rm Ran} ( \Psym ) = {\rm Ker} ( \wedge )  $ and $ {\rm Ker} ( \Psym ) = \F $.  Thus, $\F = {\rm Ran} ( 1 - \Psym ) $.  
Assumption 3. states that the corresponding  $\sigma$ is the usual flip when restricted to the center $\Ecenter$. 

We denote the restriction of the map $ \wedge $ to $\F$ by the symbol $Q$.    Thus, $Q : \F \rightarrow \twoform $ is a right $\A$-linear isomorphism. 
Let us also remark that by Proposition 6.3 of \cite{article3}, the maps $\sigma$ and $ \Psym $ are $\A$-$\A$-bilinear. For the proof of these results and more details, we refer to \cite{article3}. 
From now on, to simplify terminologies, we say that $ ( \E, d ) $ is a tame differential calculus on $\A$  if $ \E $ is the bimodule of one-forms of a tame differential calculus $ ( \Omega^\bullet(\A), d  )$.    

We list some important consequences of Definition \ref{tame} in the following proposition
(cfr. Lemma 6.2 and Theorem 3.3 of \cite{article3}):
\bppsn \label{5thdec20191} 
Suppose $ ( \E, d )  $ is a tame
differential calculus on $\A$ with $\E$ the bimodule of one-forms. 
Then the following statements hold:
\begin{itemize}
\item[1.] $\E$ is a centered bimodule, that is $ \Ecenter $ is right $\A$-total in $\E$.    
Moreover, the set $ \{ \omega \tensora \eta: \omega, \eta \in \Ecenter \} $ is both left and right $\A$-total in $\E \tensora \E$.    

\item[2.] There exists a torsionless connection $\nabla_0$ on $\E.$
\end{itemize}
\eppsn

\noindent
Let us briefly recall how  $\nabla_0$ is defined.
Since $ \E $ is finitely generated and projective as a right $ \A $-module, there exists a natural number $ n $ and an idempotent $ p \in M_n ( \A ) $ such that $ p ( \A^n ) = \E $. If $ \{ e_j: j = 1, \cdots n \} $ is a basis of the free right $ \A $-module $ \A^n, $ then the elements $ \{ \Phi_j: = p ( e_j ) : j = 1, \cdots n \} $ form a ``frame" ( in the sense of Rieffel, \cite{rieffelgh} ) of $ \E $ and $ {\rm Span}_{\mathbb{C}} \{ \Phi_j \} $ is right $\A$-total in $ \E. $
Let $ \eta $ be an element in $ \E $. Then there exists elements $ \{ a_j: j = 1, \cdots n \} $ in $ \A $ such that 
$ \eta = \sum\nolimits_j \Phi_j a_j.  $ Then the Grassmann connection $ \nabla^{Gr} $ is defined to be:
$$ \nabla^{Gr} ( \eta ) = \sum\nolimits_j \Phi_j \tensora d a_j. $$
We define $ \nabla_0: \E \rightarrow \E \tensora \E $ by the formula:
\beq\label{grass-t}
\nabla_0 = \nabla^{Gr} - Q^{-1} ( T_{\nabla^{Gr}}  ),  
\eeq
where $T_{\nabla^{Gr}}$ is the torsion of $\nabla^{Gr}$ and $Q$ is the isomorphism coming from Definition \ref{tame}.

In this set-up it is clear that $\Ecenter$ is also left $\A$-total in $\E$.    We refer to  \cite[\S 4.]{article3} for the proof  that the property of being a centered bimodules is stronger than being a central bimodule in the sense of \cite{dubois}. In particular, this means that if $\E$ is a centered bimodule,
\begin{equation} \label{16thjuly201921} 
a e = e a  \qquad  {\rm for} ~ {\rm all} ~  a \in \Acenter, ~ e \in \E. \end{equation}
Since we will be working with a tame differential calculus $ ( \E, d ) $, we are allowed to apply the properties of a centered bimodule to $\E.$

\brmrk \label{sweedler}
Suppose $ ( \E, d ) $ is a tame differential calculus on $\A$ and  let $ \xi = \sum\nolimits_i e_i \tensora f_i $ be an arbitrary element of $\E \tensora \E$.    Since $\E$ is centered ( part 1. of Proposition \ref{5thdec20191} ), there exist elements $ a_{ij} $ in $\A$ and $ h_{ij} $ in $ \Ecenter $  such that $ f_i = \sum\nolimits_j a_{ij} h_{ij} $ and hence 
$$ \xi = \sum\nolimits_{ij} e_i a_{ij} \tensora h_{ij}.$$
 Thus, if we write $ \xi $ in Sweedler's notation, $ \xi = \xi_{(0)} \tensora \xi_{(1)}$,   we can always assume that the components $ \xi_{(1)} \in \Ecenter$ without loss of generality. 
For a generic connection $\nabla$ on $\E$,  we will use the Sweedler's notation $ \nabla ( \omega ) = \omegazero \tensora \omegaone $.  However, for the very specific torsionless connection $\nabla_0$ of in Proposition \ref{5thdec20191}, we will write
$ \nabla_0 ( \omega ) = \zeroomega \tensora \oneomega $.
\ermrk

		\subsection{Pseudo-Riemannian metrics on a tame differential calculus}

\bdfn \label{6thdec20191}
Suppose $\E$ is the bimodule of one-forms of a tame differential calculus and $\sigma$ the corresponding map as defined in Definition \ref{tame}. 
  A pseudo-Riemannian metric $ g $ on $ \E $ is an element of $ {\rm Hom}_{\A} ( \E \tensora \E, \A ) $ such that
\begin{itemize} 
\item[(i)] $g$ is symmetric, that is $g \sigma = g$ ; 

\item[(ii)] the map $ \E \rightarrow \E^*, ~ \eta \mapsto V_g ( \eta ) $, with $ V_g ( \eta ) (\, \cdot \, ): \E \to \A $ defined as $ V_g ( \eta ) ( \xi ) = g ( \eta \tensora \xi ) $,  is an isomorphism of right $ \A$-modules.
\end{itemize} 
We say that a pseudo-Riemannian metric $g$ is a pseudo-Riemannian bilinear metric if in addition, $ g $ is an  $ \A $-$ \A$-bimodule map.
\edfn
 We will see that if $\E$ is the bimodule of one-forms of a tame differential calculus and $g$ is a pseudo-Riemannian bilinear metric on $\E$,  then the set $ \{ V_g ( e ) : e \in \Ecenter \} $ is a $\Acenter$-bimodule playing the role of vector fields (  Proposition \ref{prop26thjan1}   ). The following remark will be used repeatedly throughout the article.
\brmrk \label{9thdec20191}
If $g$ is a pseudo-Riemannian bilinear metric on $\E$,  the map $ V_g: \E \rightarrow \E^*  $ is 
an $ \A $-$ \A$-bimodule map. Similarly, if $\eta \in \Ecenter$, the map $ V_g ( \eta ) (\, \cdot \, ): \E \rightarrow \A $  is an $ \A $-$ \A$-bimodule map as well.
\ermrk
We have the following lemma (cf. Lemma 4.3 of \cite{article3} and Lemma 2.8 of \cite{article1}):
\blmma \label{lemma2} 
Suppose $\E$ is the bimodule of one-forms
of a tame differential calculus, $\sigma$ be  the map defined in Definition \ref{tame} and let  $g$ be a pseudo-Riemannian bilinear metric on $\E$.    
Then we have the following: 
\begin{itemize}
\item[1.] If either  $ \omega $ or $ \eta $ belongs to $ \Ecenter $, then 
$ \sigma ( \omega \tensora \eta ) = \eta \tensora \omega$. 

\item[2.]If   either of $ \omega $  or $ \eta $ belongs to $ \Ecenter $,  then    
$g ( \omega \tensora \eta ) = g ( \eta \tensora \omega )$. 

\item[3.] The element $ g ( \omega \tensora \eta ) $ belongs to  $ \mathcal{Z} ( \A ) $ if both $ \omega $ and $ \eta  $ are in 
$ \mathcal{Z} ( \E )$. 

\item[4.] If $ f $ is an element of $ \Acenter $, then $ d f \in \Ecenter $. In particular,
if $ \omega, \eta \in \Ecenter $ and $ g $ is a pseudo-Riemannian bilinear metric, then 
\begin{equation} \label{22ndjuly20192} d g ( \omega \tensora \eta ) \in \Ecenter. \end{equation}
\end{itemize}
\elmma

We next extend the metric to two-fold tensors 
 and prove an additional preparatory result which will be used in the proof of Proposition \ref{prop26thjan1}.
	\bppsn \label{vg2nondegenerate} 
 Let $ g $ be a pseudo-Riemannian bilinear metric on $\E$ where $ ( \E, d ) $ is a tame differential calculus. We define a map
\begin{align*}
& g^{(2)} :(\E\tensora \E)\tensora (\E \tensora \E) \rightarrow \A, \\
& g^{(2)}\big((\eta\tensora \xi)\tensora (\eta'\tensora \xi')\big) = g\big(\eta \tensora g(\xi \tensora \eta') \xi' \big).
 \end{align*}
 Then we have the following:
\begin{itemize}
\item[1.] The map $ V_{g^{(2)}}: \E \tensora \E \rightarrow ( \E \tensora \E )^* $   defined by
  $$ V_{g^{(2)}} ( \eta \tensora \xi ) ( \eta^\prime \tensora \xi^\prime ) = g^{(2)} \big( ( \eta \tensora \xi ) \tensora 
  ( \eta^\prime \tensora \xi^\prime )   \big) $$
 is $\A$-$\A$-bilinear and an isomorphism of right $ \A $ modules.
\item[2.] For $ \omega, \eta \in \Ecenter $, 
\begin{equation} \label{7thdec20191}  V_g ( \omega )  \tensora V_g ( \eta )  = V_{g^{(2)}}  ( \eta \tensora \omega  ).\end{equation}
\item[3.] If $\omega, \theta \in \Ecenter $ and $\xi \in \E \tensora \E$ such that $\sigma (\xi) = \xi$,  then
 \begin{equation} \label{7thdec20194} \big( V_g ( \omega ) \tensora V_g ( \theta ) \big) ( \xi ) = \big( V_g ( \theta ) \tensora V_g ( \omega ) \big) ( \xi ).\end{equation}
\end{itemize}
\eppsn
\begin{proof} The first assertion follows from Proposition 6.6 of \cite{article3} (see also Proposition 3.7 of \cite{article1}).
The equation \eqref{7thdec20191} follows by inspecting the proof of Proposition 3.7 of \cite{article1}. For  the third assertion, by Lemma 4.17 of \cite{article1}, the following equation holds for all $ x, y \in \E \tensora \E :$
 \begin{equation} \label{7thdec20195} g^{(2)} ( \Psym ( x ) \tensora y ) = g^{(2)} ( x \tensora \Psym ( y ) ).\end{equation} 
We claim that if $ \omega, \theta $ and $ \xi $ are as in \eqref{7thdec20194}, then
\begin{equation} \label{7thdec20193} g^{(2)} (  ( \omega \tensora \theta ) \tensora \xi    ) = g^{(2)} ( ( \theta \tensora \omega ) \tensora \xi ). \end{equation}
Since $ \sigma = 2 \Psym - 1 $ and \eqref{7thdec20195} holds, we get
  $$ g^{(2)} (  ( \omega \tensora \theta ) \tensora \xi    )  =  g^{(2)} (  ( \omega \tensora \theta ) \tensora \sigma \xi    ) =   g^{(2)} ( \sigma ( \omega \tensora \theta ) \tensora \xi    ) = g^{(2)} ( ( \theta \tensora \omega ) \tensora \xi ),
  $$
	which proves the claim. The equation \eqref{7thdec20194} now follows from \eqref{7thdec20193}. Indeed, using Sweedler notation $\xi = \xi_{(0)} \tensora \xi_{(1)}$ with $ \xi_{(1)} \in \Ecenter $ (Remark \ref{sweedler}), we obtain 
 \begin{align*}
  ( V_g ( \omega ) \tensora V_g ( \theta ) ) ( \xi ) & = ( V_g ( \omega ) \tensora V_g ( \theta ) ) ( \xi_{(0)} \tensora \xi_{(1)} ) = g ( \omega \tensora \xi_{(0)} ) g ( \theta \tensora \xi_{(1)} )\\
	&= g^{(2)} (  ( \theta \tensora \omega  ) \tensora ( \xi_{(0)} \tensora \xi_{(1)} ) ) =  g^{(2)} (  (  \omega \tensora \theta  ) \tensora ( \xi_{(0)} \tensora \xi_{(1)} ) ) \\ & = g ( \omega \tensora \xi_{(1)} ) g ( \theta \tensora \xi_{(0)} ) =
	g ( \theta \tensora \xi_{(0)} ) g ( \omega \tensora \xi_{(1)} ) \\ & = ( V_g ( \theta ) \tensora V_g ( \omega ) ) ( \xi ),
 \end{align*}
where we have used the fact that the component $ \xi_{(1)} \in \Ecenter $ and \eqref{7thdec20194}.
\end{proof}

\subsection{Levi-Civita connection on a tame differential calculus}

Suppose $ ( \E, d ) $ is a tame differential calculus.  We say that  a connection $\nabla$ on $\E$ is compatible with the metric $g$ on $ \Ecenter $ if for all $ \omega, \eta \in \Ecenter $, the following equation holds:
$$(g \tensora {\rm id}) \Big[ \sigma_{23}(\nabla(\omega)\tensora \eta ) + \omega \tensora \nabla ( \eta )  \Big] = d ( g ( \omega \otimes_{\Acenter} \eta ) ). 
$$
This can be extended to the whole of $ \E = \Ecenter \otimes_{\Acenter} \A  $ in the following way. 
Firstly, define a map  $ \Pi_g^0(\nabla):\Ecenter\ot\Ecenter \rightarrow \E$ by the formula: 
$$ \Pi_g^0(\nabla)(\omega \ot \eta) = (g \tensora \id ) \Big[ \sigma_{23}(\nabla(\omega)\tensora \eta ) +  
\omega \tensora \nabla ( \eta ) \Big] .
$$
It can be checked (see Section 5 of \cite{article3}) that $ \Pi_g^0(\nabla) $ descends to a map on $\E \tensora \E$. 
Thus the compatibility of  $ \nabla $ is with $ g $ on $ \Ecenter $  can be written as 
$$ 
\Pi_g^0 ( \nabla ) ( \omega \otimes_{\Acenter} \eta ) = d ( g ( \omega \tensora \eta ) ) \qquad \forall ~ \omega, \eta \in \Ecenter. 
$$
In Subsection 4.1 of \cite{article1}, it was proven that the when $ \E = \Ecenter \otimes_{\Acenter} \A  $ one can define a canonical extension $ \Pi_g ( \nabla ) : \E \tensora \E \rightarrow \E $ of the map $ \Pi^0_g ( \nabla ) $. Concretely, for all $ \omega, \eta \in \Ecenter $ and $ a \in \A $,  
one has
\begin{equation} \label{22ndjuly20195} \Pi_g ( \nabla ) ( \omega \tensora \eta a ) = \Pi^0_g ( \nabla )  ( \omega \otimes_{\Acenter} \eta ) a + g ( \omega \tensora \eta ) da. \end{equation}
 We say that a connection $ \nabla $ is compatible with $ g $ on the whole of $ \E $ if for all $\eta, \xi$ in $\E$, it holds that 
 \begin{equation} \label{10thaugust20191} \Pi_g ( \nabla ) ( \eta \tensora \xi ) = d ( g ( \eta \tensora \xi ) ). \end{equation} 

Finally, the next theorem is the main result of \cite{article3}.
\bthm \label{existenceuniqueness}
Suppose $ ( \Omega^\bullet ( \A ), d ) $ is a tame differential calculus as in Definition \ref{tame} on  $ \A $.  
If $ g $ is a pseudo-Riemannian bilinear metric on $ \E $, then there exists a unique connection $\nabla$ on $ \E $ which is torsionless and compatible with $ g. $
\ethm

\section{A class of  derivations from a tame differential calculus} \label{section3}
	
In this section we show that if $ ( \E, d ) $ is a tame differential calculus as in Definition \ref{tame}, then  there is a right $\A$-total $ \Acenter $-submodule 
of $ \E^* $ which can be viewed as the analogue of vector fields in this context. We will denote this submodule by the symbol $ \mathcal{X} ( \A )$. The elements of $ \mathcal{X} ( \A ) $  act by derivations on $ \A $. However, our vector fields will be acting on the differential forms and not the other way around. The $ \Acenter $-bimodule of vector fields $ \mathcal{X} ( \A ) $ is defined via the  following: 
	\blmma \label{lemma26thjan1}
	Let $ ( \E, d ) $ be a tame differential calculus and let
	$ g $ be a pseudo-Riemannian bilinear metric on $ \E $. Define,
	\beq\label{6thdec20192}
	\mathcal{X} ( \A ) =  \{ V_g ( \omega ) : \omega \in \Ecenter \} \subseteq \E^* \, .
	\eeq 
	Then, 
	 \begin{itemize}
	 \item[1.]
	 If $ X \in \mathcal{X} ( \A ) $ and $ \eta \in \Ecenter $, then $ X ( \eta ) $ belongs to $\Acenter. $
	
	\item[2.] $ \mathcal{X} ( \A ) = \mathcal{Z} ( \E^* ). $
	
	\item[3.] $ \mathcal{X} ( \A )  $ is a right $\Acenter$ submodule of $\E^*$ which is right $\A$-total in $\E^*.$
	
\item[4.] The map $ \delta_\phi $ defined for an element  $ \phi \in \E^* $ by
	$$ \delta_\phi : \A \rightarrow \A, ~ \delta_\phi ( a ) = \phi ( da ). $$
	 acts by derivations on $\A$  if and only if $ \phi \in \mathcal{X} ( \A ). $
	\end{itemize}
	\elmma
	\begin{proof} For part 1.: let $ X = V_g ( \omega ) $ for some $ \omega \in \Ecenter$,  then $ X ( \eta ) = g ( \omega \tensora \eta )  $ belongs to $ \Acenter $ by part 3. of Lemma \ref{lemma2}. For part 2.: use the $\A$-$\A$-bimodule structure of $ \E^* = \Hom_\A ( \E, \A ) $ as spelled out in \eqref{7thdec20196}. The fact that $ \mathcal{X} ( \A ) \subseteq \mathcal{Z} ( \E^* ) $ is a simple consequence of the bilinearity of $g$.    Conversely, suppose $\phi$ is an element of $\E^*$.    Since $ V_g: \E \rightarrow \E^* $ is an isomorphism, there exists $\omega $ in $\E$ such that $\phi = V_g ( \omega )$.    As $ \phi \in \mathcal{Z} ( \E^* )$,  then for all $ a \in \A $, we must have 
		$$ a V_g ( \omega ) = V_g ( \omega ) a. $$
		However, $ a V_g ( \omega ) = V_g ( a \omega ) $ by the bilinearity of $g$ and therefore, $ V_g ( a \omega ) = V_g ( \omega a ) $.  This implies that $ a \omega = \omega a $ as $V_g$ is invertible and so $ V_g ( \omega ) \in \mathcal{X} ( \A ) $.  This proves point 2.
		For the third assertion: since the differential calculus is tame, Proposition \ref{5thdec20191} implies that $\Ecenter \A = \E $.  Since $V_g$ is a right $\A$-linear isomorphism, we get $ \mathcal{X} ( \A ) \A = \E^*$.     
		
		Finally, we prove the fourth assertion. If $\phi $ belongs to $\mathcal{X} ( \A ), $ then the Leibniz rule for the differential $d$ easily implies that $\delta_\phi$ is a derivation.
Conversely, if $ \phi \in \E^* $ is such that $ \delta_\phi $ is a derivation, then for all $a, b \in \A $, we have
		\begin{align*} ( a \phi ) ( db ) & = a \delta_\phi ( b ) = \delta_\phi ( a b ) - \delta_\phi ( a ) b = \phi ( d ( a b ) ) - \phi ( da ) b \\ & = \phi ( d ( a b ) - da b ) = \phi ( a db ) = ( \phi a ) ( db ).
		\end{align*}
		Therefore, for all $ a, b \in \A, $
		$$ ( a \phi ) ( db ) = ( \phi a ) ( db ).  $$
		Since $ a \phi $ and $ \phi a $ are right $\A$-linear, this implies that $ \phi \in \mathcal{Z} ( \E^* ) = \mathcal{X} ( \A ).$
		\end{proof}
		
		Next proposition tells us that the right $\Acenter$-submodule $ \mathcal{X} ( \A ) $ is naturally a Lie subalgebra of the set of all derivations from $\A$ to $\A$
 and so it can be viewed as an analogue of the Lie algebra of vector fields. 

	\bppsn \label{prop26thjan1} 
Let $ ( \E, d ) $ be a tame differential calculus.
Given any $X, Y \in \mathcal{X} ( \A )$,  there exists a unique element $ [ X, Y ] \in \mathcal{X} ( \A ) $ such that
\begin{equation} \label{26thjan1} \delta_{[X, Y]} = \delta_X \circ \delta_Y - \delta_Y \circ \delta_X \, . 
\end{equation}
Thus $ \mathcal{X} ( \A ) $ is a Lie algebra under the Lie  bracket defined by the previous equation.   
\eppsn
  \begin{proof} 
		
	\noindent
	For proving \eqref{26thjan1}, let us first note that since the set of all derivations is closed under Lie-bracket, the expression $  [ \delta_X, \delta_Y ] $ makes sense and is again a derivation. We need to show that this derivation is of the form $ \delta_Z $ for some $ Z \in \mathcal{X} ( \A ) $.  To this end, we define 
	$$ Z: \E \rightarrow \A \qquad {\rm by} \qquad Z ( \sum\nolimits_i d a_i b_i ) = \sum\nolimits_i [ \delta_X, \delta_Y ] ( a_i ) b_i. $$
	If $ Z $ is well defined, it is immediate that $ Z $ is right $\A$-linear, that is $ Z \in \E^* $.  Moreover, it is also clear that $ \delta_Z = [ \delta_X, \delta_Y ]$.    Now, since $   [ \delta_X, \delta_Y ] $ is a derivation, the fourth assertion of Lemma \ref{lemma26thjan1} implies that  $ Z \in \mathcal{X} ( \A ) $.    So we are left with proving that $ Z $ is well-defined.
	
	\noindent 
	Suppose $ a_i, b_i \in \A $ be such that $ \sum\nolimits_i d a_i b_i = 0 $.  Then $ X ( \sum\nolimits_i d a_i b_i ) = Y ( \sum\nolimits_i d a_i b_i ) = 0 $ and hence
	\begin{equation} \label{31stjan1} \sum\nolimits_i X ( d a_i ) b_i = \sum\nolimits_i Y ( d a_i ) b_i = 0. \end{equation} 
	Moreover,  
	$$ 0 = d ( \sum\nolimits_i d a_i b_i ) = \sum\nolimits_i d a_i \wedge d b_i = \wedge ( \sum\nolimits_i d a_i \tensora d b_i ),  $$
	that is,  $ \sum\nolimits_i d a_i \tensora d b_i \in {\rm Ker} ( \wedge  ) = {\rm Ran} ( \Psym ) $ (see Definition \ref{tame}) and so 
	$$ \sigma ( \sum\nolimits_i d a_i \tensora d b_i ) = ( 2 \Psym - 1 ) ( \sum\nolimits_i d a_i \tensora d b_i )   = \sum\nolimits_i d a_i \tensora d b_i. $$ 
	So if $ X = V_g ( \omega ) $ and $Y = V_g ( \eta ) $ for some $ \omega, \eta \in \Ecenter$,  then \eqref{7thdec20194} implies that 
	\begin{multline} \label{31stjan2}
	( X \tensora Y - Y \tensora X ) ( \sum\nolimits_i d a_i \tensora d b_i ) \\
	 = ( V_g ( \omega ) \tensora V_g ( \eta )   ) ( \sum\nolimits_i d a_i \tensora db_i )- ( V_g ( \eta ) \tensora V_g ( \omega )   ) ( \sum\nolimits_i d a_i \tensora db_i )\\            
	=   ( V_g ( \eta ) \tensora V_g ( \omega )   ) ( \sum\nolimits_i d a_i \tensora db_i ) - ( V_g ( \eta ) \tensora V_g ( \omega )   ) ( \sum\nolimits_i d a_i \tensora db_i ) = 0.
	\end{multline}
	Now,
	\begin{align*}
	 \sum\nolimits_i [ \delta_X, & \delta_Y ] ( a_i ) b_i =  \sum\nolimits_i [ \delta_X ( Y ( d a_i ) ) b_i - \delta_Y ( X ( d a_i ) ) b_i ]\\
	&=  \delta_X ( \sum\nolimits_i Y ( d a_i ) b_i ) - \sum\nolimits_i Y ( d a_i ) \delta_X ( b_i )   -  \delta_Y ( \sum\nolimits_i X ( d a_i ) b_i ) + \sum\nolimits_i X ( d a_i ) \delta_Y ( b_i )\\
	&= ( X \tensora Y - Y \tensora X ) ( \sum\nolimits_i d a_i \tensora d b_i ) ~ {\rm ( } ~ {\rm by} ~ \eqref{31stjan1} ~ {\rm ) }
	\end{align*}
	which is $ 0 $ by \eqref{31stjan2}. This shows  that $ Z $ is well-defined, hence completes the proof of the proposition.   \end{proof}

Thus for a tame differential calculus $ ( \E, d ), $ we have defined a Lie subalgebra $ \mathcal{X} ( \A ) $ 
of the Lie algebra $ {\rm Der} ( \A )$. Before we proceed, we relate our approach to the one of derivation based differential calculus in the spirit of that proposed in \cite{dubois88}. We shall compute this $\Acenter$-submodule explicitly for two classes  of differential calculi. 
\bxmpl
 Let $\A$ be an algebra and let $L$ be a Lie algebra which acts on $ A $ by derivations. If the Chevalley
 differential calculus $ ( C^* ( L, \A ), d ) $ is tame, then $ \mathcal{X} ( \A ) \simeq \Acenter \tensorc L$.  

By definition the differential calculus $ ( C^* ( L, \A ), d ) $ is the graded differential algebra of $L$ with coefficients in $\A$. 
In particular $ C^0 ( L, \A ) = \A $ and the space 
$ C^1 ( L, \A ) $ of one-forms is the space of all $\IC$-linear maps 
from $L$ to $A$ with its $\A$-$\A$-bimodule structure defined by:
$$ 
( \phi a  ) ( X ) = \phi ( X ) a \, , \qquad ( a \phi ) ( X ) = a \phi ( X ) 
$$
for all $a \in \A$, $\phi \in C^1 ( L, \A )$ and  $X \in L$.
The differential $d: \A \rightarrow C^1 ( L, \A ) $ is $da (X) = X(a)$. 

As per our convention, we denote by $\E = C^1 ( L, \A ) $ the space of one-forms. 
The fact that the calculus is tame means in particular that 
$\E = \Ecenter \otimes_{\Acenter} \A$.
Thus, any element of $ \E $ is a unique right $\A$-linear combination of elements 
$\theta_j $ in $\E,$ where $\theta_j $ is defined as 
$$ 
\theta_j  ( X_k ) = \delta_{j k} \, 1_A ,
$$
for a given basis $ \{  X_k, \, k=1, \cdots, \dim(L)\} $ of $L$. Therefore, $\E$ is a free right $\A$-module with basis $ \{ \theta_j \}$. Moreover, the $\A$-$\A$- bimodule structure on $\E$ also implies that $ \theta_j $ belongs to $ \Ecenter$. 
Next, let us introduce elements $ \xi_k \in \E^* $ defined by 
$$ 
\xi_k \, ( \sum\nolimits_j \theta_j a_j ) = a_k.  
$$
 It follows that $ \E^* $ is a free right $\A$-module with basis $ \{ \xi_k: k = 1, \cdots, \dim(L)\}$.  
Now, since we are assuming that the differential calculus $ ( C^* ( L, \A ), d ) $ is tame, then part 4. of Lemma \ref{lemma26thjan1} implies that $ \mathcal{X} ( \A ) $ is equal to $ \mathcal{Z} ( \E^* )$. 
It can be easily checked that an element $ \sum\nolimits_i \xi_j a_j $ belongs to $ \mathcal{Z} ( \E^* ) $ 
if and only if $ a_j $ belongs to $ \Acenter $ for all $j$. Hence, 
$$ 
\mathcal{X} ( \A ) = \mathcal{Z} ( \E^* ) = \Acenter \tensorc L. 
$$
As a final point we remark that the isomorphism $ \E \rightarrow \E^*$ as in point (ii) of Definition \ref{6thdec20191} 
can be defined via the Euclidean metric given on generators by
$g (\theta_k \tensora \theta_j) =\delta_{k j}$; the isomorphism is then just given by $ \theta_k \mapsto \xi_k$, 
for $k = 1, \cdots, \dim(L)$.  
\exmpl

\bxmpl
If  the derivation based differential calculus $ ( \Omega^*_{{\rm Der}} ( \A ), d  ) $ as given in \cite{dubois88} is tame, then $ \mathcal{X} ( \A ) $ is isomorphic to the set of all derivations $ {\rm Der} ( \A ). $
 
Given the algebra $\A$, the pair $(\Omega^1_{{\rm univ}} ( \A ), d_{{\rm univ}})$ denotes its universal 
differential calculus. Its universal property implies that given a derivation 
 $ X \in {\rm Der} ( \A )$ there exists a unique $\A$-bimodule homomorphism $i_X: \Omega^1_{{\rm univ}} (\A) \rightarrow \A$ 
 such that $i_X \circ d_{{\rm univ}} = X$. 
Next, consider the ideal $ I = \{ \omega \in \Omega^1_{{\rm univ}} ( \A ): i_X (\omega) 
= 0 \, \, \mbox{for all} \, \, X \in {\rm Der} ( \A )   \} $ and the quotient 
$$
\dubois = \Omega^1_{{\rm univ}} ( \A )  / I 
$$
and denote by $d$ the differential $d_{{\rm univ}}$ when restricted to $\dubois$. The differential extends 
canonically as a derivation $d: \Omega^k_{{\rm Der}}(\A) \rightarrow \Omega^{k + 1}_{{\rm Der}}(\A)$, where 
$$
\Omega^k_{{\rm Der}}(\A) = \dubois \otimes_{\A} \dubois \cdots \otimes_{\A} \dubois \quad (k \,\, \mbox{factors}) 
$$
with $d^2=0$, thus giving rise to a differential calculus.
There is a residual universality for $\dubois$ in that given a derivation 
 $ X \in {\rm Der} ( \A )$ there exists a unique $\A$-bimodule 
 map $\phi_X$ from $ \dubois $ to $\A$ such that $ X ( a ) = \phi_X ( d ( a ) )$.
 Now, in the hypothesis that the calculus is tame, with $\E = \dubois$, from $\E = \Ecenter \otimes_{\Acenter} \A$ 
one sees that the set of all $\A$-$\A$-bilinear maps from $\E $ to $\A$ is equal to $ \mathcal{Z} ( \E^* )$.
Also, by part 2. of Lemma \ref{lemma26thjan1}, $ \mathcal{X} ( \A ) = \mathcal{Z} ( \E^* )$. It follows that the map $ X \mapsto \phi_X  $ is a one-one map from $ {\rm Der} ( \A ) $ to $ \mathcal{X} ( \A )$. This map is also onto: if $ \phi \in \mathcal{X} ( \A ) = \mathcal{Z} ( \E^* ),  $ then by part 4. of Lemma \ref{lemma26thjan1},
$$ 
X_\phi: \A \rightarrow \A, \quad X_{\phi} ( a ) = \phi ( d ( a ) ) 
$$
is a derivation on $\A.$ Hence, $ \mathcal{X} ( \A ) $ is isomorphic to $ {\rm Der} ( \A )$.
\exmpl 
	
Interestingly, in our context one can make sense of the Lie-bracket  $ [ X, \phi ] $ where $ X  $ belongs to  $ \mathcal{X} ( \A ) $ and $ \phi $ is a general element of  $ \E^* $.  We will actually need this fact for the statement of Theorem \ref{koszulformula2} below. To that end, let us first observe that since $ ( \E, d ) $ is tame, we can use the identification $ \E = \Ecenter \otimes_{\Acenter} \A $ and right $\A$-linearity of the map $ V_g $ to write
\begin{align*}
 \E^* &=V_g ( \E ) = ( V_g \otimes_{\Acenter}  {\rm id}  ) ( \Ecenter \otimes_{\Acenter} \A ) \\
      & = V_g ( \Ecenter ) \otimes_{\Acenter} \A = \mathcal{X} ( \A ) \otimes_{\Acenter} \A \\
      & = \mathcal{Z} ( \E^* ) \otimes_{\Acenter} \A             
\end{align*}
and we have used the second assertion of Lemma \ref{lemma26thjan1}. The identification $ \E^* = \mathcal{Z} ( \E^* ) \otimes_{\Acenter} \A $ will be used for the next definition.
 \bdfn
Suppose $ X \in \mathcal{X} ( \A ) $ and $  \phi = \sum\nolimits_i Y_i \otimes_{\Acenter} a_i  \in \E^* $ for some $ Y_i \in \mathcal{X} ( \A ) $ and $ a_i \in \A $.  We define
\begin{equation} \label{4thfeb20192} [ X, \phi ] := \sum\nolimits_i ( [ X, Y_i ] a_i + \delta_X ( a_i ) Y_i ). \end{equation}
\edfn
\noindent
We need to prove that the above \eqref{4thfeb20192} is well-defined.
\bppsn \label{4thfeb2019prop}
Suppose $ X $ and $ \phi $ be as above. Then $ [ X, \phi ] $ is well-defined.
\eppsn
\begin{proof} Suppose $ \sum\nolimits_i Y_i \otimes_{\Acenter} a_i = 0 $.  We need to show that $  \sum\nolimits_i ( [ X, Y_i ] a_i + \delta_X ( a_i ) Y_i ) = 0.$
   We claim that it is enough to prove that for all $ b \in \A, $
	\begin{equation} \label{4thfeb2019} \sum\nolimits_i ( [ X, Y_i ] a_i + \delta_X ( a_i ) Y_i ) ( db ) = 0.  \end{equation}
	Indeed, if \eqref{4thfeb2019} holds, then for all $ c \in \A $, 
	$$ \sum\nolimits_i \big( [ X, Y_i ] a_i + \delta_X ( a_i ) Y_i \big) ( (db ) c )) = \sum\nolimits_i \big( [ X, Y_i ] a_i ( db ) + \sum\nolimits_i \delta_X ( a_i ) Y_i  ( db ) \big) c = 0$$
	by the right $ \A $-linearity of both $ [ X_i, Y_i ] $ and $ Y_i. $
	As for \eqref{4thfeb2019}: since $ [ X, Y_i ] \in \mathcal{X} ( \A ) = \mathcal{Z} ( \E^* ) $ by part 2. of Lemma \ref{lemma26thjan1}, we obtain 
	\begin{align*}
  	\sum\nolimits_i  [ X, Y_i ] a_i ( db )&=  \sum\nolimits_i a_i [ X, Y_i ] ( db ) = \sum\nolimits_i a_i [ \delta_X, \delta_{Y_i} ] ( b ) ~ {\rm ( } ~ {\rm by} ~ \eqref{26thjan1} ~    {\rm ) }\\
	&= \sum\nolimits_i a_i [ \delta_X ( \delta_{Y_i} ( b ) ) - \delta_{Y_i} (  \delta_X ( b )  )  ]\\
	&= \sum\nolimits_i a_i \delta_X ( Y_i ( d b ) ) - \sum\nolimits_i a_i Y_i ( d ( \delta_X ( b ) ) )\\
	&= \sum\nolimits_i a_i \delta_X ( Y_i ( d b ) ) ~ {\rm ( } ~ {\rm as} ~ \sum\nolimits_i a_i Y_i = \sum\nolimits_i Y_i a_i = 0 ~    {\rm ) }
	\end{align*}
	and hence
	\begin{align*} \sum\nolimits_i ( [ X, Y_i ] a_i + \delta_X ( a_i ) Y_i ) ( db ) & = \sum\nolimits_i a_i \delta_X ( Y_i ( d b ) ) + \sum\nolimits_i \delta_X ( a_i ) Y_i ( db ) \\ & = \delta_X ( ( \sum\nolimits_i a_i  Y_i ) ( d b ) )   = 0.
	\end{align*}
This proves \eqref{4thfeb2019}.
\end{proof}

	For the rest of the subsection, we discuss some consequences of the definitions and observations above. The  next lemma is needed in the proof of Proposition \ref{prop26thjan2}.
	
	\blmma \label{lemma26thjan4}
	If $ \nabla_0 $ is the torsionless connection of Proposition \ref{5thdec20191} as given in \eqref{grass-t}, and $X, Y \in \mathcal{X} ( \A ) $, then for all $a \in \A$,  the following equation holds:
	$$ ( X \tensora Y - Y \tensora X ) \nabla_0 ( da ) = 0. $$
	\elmma
	\begin{proof} Since $ {\rm Ran} ( \Psym ) = {\rm Ker} ( \wedge ) $ (see Definition \ref{tame}), we have
		  $$ \wedge \nabla_0 ( \omega ) = \wedge ( \Psym \nabla_0 ( \omega ) + (  1 - \Psym ) \nabla_0 ( \omega ) )  = \wedge (  1 - \Psym ) \nabla_0 ( \omega ) \qquad  {\rm for} ~ {\rm all} ~ \omega ~ \in ~ \E. $$
			Since $ \nabla_0 $ is torsionless, $ \wedge \nabla_0 ( d a ) = - d ( d a ) = 0 $.  
			Thus, $ \wedge ( 1 - \Psym ) \nabla_0 ( da ) = 0 $ and, since the map $ \wedge $ is an isomorphism from $ {\rm Ran} ( 1 - \Psym ) = \F $ onto $ \Omega^2 ( \A ) $, we can conclude that
			$$ 0 = ( 1 - \Psym ) \nabla_0 ( da ) = \frac{( 1 - \sigma )}{2} \nabla_0 ( da ). $$
		If we write $ \nabla_0 ( da ) = \zeroomega \tensora \oneomega $ (with $ \oneomega \in \Ecenter $ as in Remark \ref{sweedler}), the above equation implies that
			$ \zeroomega \tensora \oneomega = \oneomega \tensora \zeroomega. $
		Thus, 
		 \begin{align*}
		   ( X \tensora Y - Y \tensora X ) \nabla_0 ( da ) &= ( X \tensora Y )  ( \zeroomega \tensora \oneomega ) - ( Y \tensora X ) ( \oneomega \tensora \zeroomega )\\ 
			&=  X ( \zeroomega ) Y ( \oneomega ) - Y ( \oneomega ) X ( \zeroomega )\\
			&= X ( \zeroomega ) Y ( \oneomega ) - X ( \zeroomega ) Y ( \oneomega )\\
			&= 0
		 \end{align*}
			and we have used part 1. of Lemma \ref{lemma26thjan1} to observe that $ Y ( \oneomega ) \in \Acenter. $
		 \end{proof}

	 For the proof of Proposition \ref{torsionless26thjan}, it will be helpful to have the following classical formula involving $\nabla_0$,  for the Lie bracket $ [ X, Y ] $ of two elements $X, Y$ of $\mathcal{X} ( \A )$.
	\bppsn \label{prop26thjan2}
	Suppose $ ( \E, d ) $ is a tame differential calculus. Let   $X, Y$ be elements of $\mathcal{X} ( \A )$ and $ \xi $ be an element of $ \E $.  If $ \nabla_0 $ is the torsionless connection of Proposition \ref{5thdec20191}, then
	\begin{equation} \label{26thjan9} [ X, Y ] ( \xi ) = X ( d ( Y ( \xi ) ) ) - Y ( d ( X ( \xi ) ) ) + ( X \tensora Y - Y \tensora X ) \nabla_0 ( \xi ).\end{equation}
	\eppsn
	\begin{proof}  Let us define a map $ \Psi ( X, Y ) : \E \rightarrow \A $ by the right hand side of \eqref{26thjan9}.
	We claim that this map $ \Psi ( X, Y )  $ is right $\A$-linear. Indeed, for elements $a, b \in \A$,  we compute
	\begin{equation} \label{9thaugust20191} \Psi ( X, Y ) ( da ) = X ( d ( Y ( da ) )  ) - Y ( d ( X ( da ) ) ) + 0 \end{equation}
	by an application of Lemma \ref{lemma26thjan4}. Moreover, 
	\begin{align*}
	 \Psi ( X, Y ) ( da b ) &= X ( d ( Y ( da b ) ) ) - Y ( d ( X ( da b ) ) ) + ( X \tensora Y - Y \tensora X ) \nabla_0 ( da b )\\
	&= X ( d ( Y ( da ) ) b ) + X ( Y ( da ) db ) - Y ( d ( X ( da ) ) b ) - Y ( X ( da ) db )\\
		& \quad + ( X \tensora Y ) ( \nabla_0 ( da ) b + da \tensora db ) - ( Y \tensora X ) ( \nabla_0 ( da ) b + da \tensora db )\\
	&= [ X ( d ( Y ( da ) ) ) - Y ( d ( X ( da ) ) ) + ( X \tensora Y - Y \tensora X ) \nabla_0 ( da )] b\\
	&= \Psi ( X, Y ) ( da ) b. 
	\end{align*}
	Here, we have used Lemma  \ref{lemma26thjan4} and that the elements $ X, Y $ belonging to $ \mathcal{X} ( \A ) = V_g ( \Ecenter ) $ are both left and right $\A$-linear by virtue of the second assertion in Remark \ref{9thdec20191}.
	 Since $ \E = {\rm Span} \{ da b : a, b \in \A \} $, this proves the claim. 
	However, \eqref{9thaugust20191} implies that 
	$$ \Psi ( X, Y )  ( da ) = \delta_X ( \delta_Y ( a ) ) - \delta_Y ( \delta_X ( a ) ) = [ \delta_X, \delta_Y ] ( a ) = \delta_{[ X, Y ]} ( a ) = [ X, Y ] ( da ).  $$
	Since $ [ X, Y ] $ is an element of $\E^*$,  it is right $\A$-linear while $ \Psi ( X, Y ) $ is right $ \A $-linear by our claim. Thus, the equation \eqref{26thjan9} holds for all $ \xi $ in $\E$.     \end{proof}

\section{The covariant derivative from a connection} \label{section4}

The analysis made in the previous section allows us to define a covariant derivative on $ \mathcal{X} ( \A ) $ starting from a connection on one-forms. Moreover, the connection on one-forms can be recovered from the covariant derivative on $ \mathcal{X} ( \A ) $.   
	
\bdfn
	Suppose $( \E, d )$ is a tame differential calculus and $ \nabla $ is a connection on $ \E $.  Given elements  $ X, Y $ in $ \mathcal{X} ( \A ) $,  we define  $ \nabla_Y X \in \E^* $ by the equation
	\begin{equation} \label{30thjan3} ( \nabla_Y X ) ( \omega ) = \delta_Y ( X ( \omega ) ) - ( X \tensora Y ) ( \nabla ( \omega ) ). \end{equation}
	\edfn
	\noindent It can be easily seen that $ \nabla_Y X $ is indeed an element of $ \E^* $, that is 
	$\nabla_Y X ( \omega a ) = ( \nabla_Y X ( \omega ) ) a  $ for all $ X, Y \in \mathcal{X} ( \A ), \, 
	\omega \in \E $ and $ a \in \A. $
	
	\blmma \label{9thnov20192}
	Suppose  $ \nabla $ is a connection on $\E$ where $ ( \E, d ) $ is a tame differential calculus. If $X, Y, X^\prime, Y^\prime \in \mathcal{X} ( \A )$ and $ a \in \Acenter $, then the following equations hold:
	\begin{align} \label{13thapril20191} \nabla_Y ( X + X^\prime ) & = \nabla_Y X + \nabla_Y X^\prime, \qquad \nabla_{Y + Y^\prime} X = \nabla_Y X + \nabla_{Y^\prime} X, \nonumber \\
	 \nabla_{Y a} X &= ( \nabla_Y X ) a, \qquad \nabla_{Y} ( X a ) = ( \nabla_ Y X ) a + X \delta_Y ( a ). \end{align}
	\elmma
	\begin{proof} The first and the second equalities are straightforward. For proving the fourth equality, we see that if $ \omega $ belongs to $\E$,  then
	\begin{equation} \label{10thdec20191} \nabla_Y ( X a ) ( \omega ) = \delta_Y ( X a ( \omega ) ) - ( X a \tensora Y ) (  \nabla ( \omega ) ). \end{equation}
	Now, using Sweedler notation to write $ \nabla ( \omega ) = \omegazero \tensora \omegaone $, we observe that
	\begin{align*}
	( X a \tensora Y ) (  \nabla ( \omega ) ) &= X a ( \omegazero ) Y ( \omegaone ) = X ( a \omegazero ) Y ( \omegaone ) ~ {\rm (} ~ {\rm by} ~ \eqref{7thdec20196} ~ {\rm )}\\
	                                          &= X ( \omegazero a ) Y ( \omegaone ) =  X ( \omegazero  ) a Y ( \omegaone )
 \\                                              
	                                          & \qquad\qquad\qquad{\rm (} ~ {\rm by} ~ \eqref{16thjuly201921}~ {\rm and} ~ {\rm the} ~ {\rm right} ~ \A- {\rm linearity} ~ {\rm of} ~ X ~ {\rm )}\\                                              
	                                          &= X ( \omegazero  )  Y ( \omegaone ) a = ( X  \tensora Y ) (  \nabla ( \omega ) ) a 
	                                 \end{align*} 
  as $a  \in \Acenter.$	Therefore, from \eqref{10thdec20191}, we obtain
	\begin{align*}
 	\nabla_Y & ( X a ) ( \omega ) = \delta_Y ( X a ( \omega ) ) - ( X  \tensora Y ) (  \nabla ( \omega ) ) a \\
	                            &= \delta_Y ( X a ( \omega ) ) - ( X  \tensora Y ) (  \nabla ( \omega ) ) a - X ( \omega ) Y ( da ) + X ( \omega ) Y ( da ) \\
&= \delta_Y ( X a ( \omega ) ) - ( X  \tensora Y ) (  \nabla ( \omega ) ) a - ( X \tensora Y ) ( \omega \tensora da ) + X ( Y ( da ) \omega )\\
& \qquad\qquad\qquad {\rm (} ~ {\rm since} ~ X ~ {\rm is} ~ {\rm left} ~ \A-{\rm linear} ~ {\rm )}\\
&= \delta_Y ( X a ( \omega ) ) - ( X  \tensora Y ) (  \nabla ( \omega a ) ) + X ( Y ( da ) \omega )\quad  {\rm (} ~ {\rm applying} ~ {\rm the}~ {\rm Leibniz} ~ {\rm rule} ~ {\rm for} ~ \nabla ~ {\rm )}\\
&= \delta_Y ( X a ( \omega ) ) - ( X  \tensora Y ) (  \nabla ( a \omega  ) ) + X ( Y ( da ) \omega ) ~ {\rm (} ~ {\rm by} ~ \eqref{16thjuly201921}  ~ {\rm )}\\
															&= ( \nabla_Y X  ) ( a \omega ) + X (  \delta_Y ( a ) \omega )\\
															&= ( ( \nabla_ Y X ) a + X \delta_Y ( a ) ) ( \omega ).
	                 \end{align*}
									Therefore, for all $ \omega $ in $\E$,  we get
									$$ \nabla_{Y} ( X a ) ( \omega ) = ( ( \nabla_ Y X ) a + X \delta_Y ( a ) ) ( \omega ) $$
									which proves the fourth equality.
	Finally, the third equality follows easily by applying \eqref{16thjuly201921}.  \end{proof}
	
	\bppsn \label{15thapril2019prop2}
Suppose $( \E, d )$ is a tame differential calculus and $ \nabla_1,  \nabla_2 $ are two connections on $ \E $ such that for all $ X, Y \in \mathcal{X} ( \A )$,  it holds that $( \nabla_1 )_Y X = ( \nabla_2 )_Y X $.  Then $ \nabla_1 = \nabla_2 $.  
	\eppsn
	\begin{proof} The equation \eqref{30thjan3} implies that for all $ X, Y \in \mathcal{X} ( \A ) $ and for all $ \omega \in \E, $
	 $$ ( X \tensora Y ) ( \nabla_1 ( \omega ) ) = ( X \tensora Y ) (  \nabla_2 ( \omega )  ). $$
	Therefore, for all $ \eta, \theta \in \Ecenter $, the following equation holds:
	 $$ ( V_g ( \theta ) \tensora V_g ( \eta ) ) (  (  \nabla_1 - \nabla_2 ) ( \omega ) ) = 0. $$
	Hence, \eqref{7thdec20191} implies that 
	$$ \vgtwo ( \eta \tensora \theta ) (  (  \nabla_1 - \nabla_2 ) ( \omega )  ) = 0.$$
By part 1. of Proposition \ref{vg2nondegenerate}, the map $V_{g^{(2)}} : \E \tensora \E \rightarrow ( \E \tensora \E )^* $ is an isomorphism of right $\A$-modules which is also left $\A$-linear.  As $ \{ \eta \tensora \theta: \eta, \theta \in \Ecenter \}$ is left $\A$-total in $ \E \tensora \E $ (part 1. of Proposition \ref{5thdec20191}), we can conclude that $ (  \nabla_1 - \nabla_2 ) ( \omega ) = 0.$
		\end{proof}
	
Classically, given the covariant derivative $\nabla_Y X$, for $X, Y \in \mathcal{X} ( \A )$, one reconstructs the connection using \eqref{30thjan3}. In this noncommutative set-up, one recovers the connection provided $ \nabla_Y X $ belongs to $ \mathcal{X} ( \A ) $ for all $ X, Y \in \mathcal{X} ( \A ) $.  This is the content of the next proposition.

   \bppsn \label{15thapril2019prop}
 Suppose we are given a collection $ \{  \nabla_Y X: X, Y \in \mathcal{X} ( \A )   \} \subseteq \mathcal{X} ( \A ) $ such that equations \eqref{13thapril20191}  
 are satisfied.
Then there exists a unique connection $ \nabla $ on $ \E $ such that for all $ X, Y \in \mathcal{X} ( \A ) $ and for all $ \omega \in \E $, 
 $$ ( X \tensora Y ) ( \nabla ( \omega ) ) = \delta_Y ( X ( \omega ) ) - ( \nabla_Y X ) ( \omega ). $$
	\eppsn
\noindent 
Our  goal is to prove Proposition \ref{15thapril2019prop}. We will need some auxiliary results. Throughout the rest of the  section, we work under the hypotheses of Proposition \ref{15thapril2019prop}.
	
	\blmma \label{13thapril2019lemma1}
	Let $ \omega \in \E $ and $ \theta, \eta \in \Ecenter $.  Define an element $ \widetilde{T_{\eta, \theta}} ( \omega ) \in \A $ by the formula:
	\begin{equation} \label{13thapril20193} \widetilde{T_{\eta, \theta}} ( \omega ) = \delta_{V_g ( \eta )} ( g ( \theta \tensora \omega )  ) - \nabla_{V_g ( \eta )} V_g ( \theta ) ( \omega ). \end{equation}
	Then for all $ a \in \Acenter$, it holds that $\widetilde{T_{\eta, \theta a}} ( \omega  ) = \widetilde{T_{\eta, \theta}} ( \omega ) a $.  
	\elmma
	\begin{proof} Since $ \omega $ belongs to $\E$ and $ \theta $ belongs to $\Ecenter$,  point 2. of Lemma \ref{lemma2} implies $ g (  \theta \tensora \omega ) = g ( \omega \tensora \theta )$.    Moreover, as $ a $ belongs to $\Acenter$,  the element $da$ belongs to $ \Ecenter $ (part 4. of Lemma \ref{lemma2}) and therefore, part 3. of Lemma \ref{lemma2} implies that $ g ( \eta \tensora da ) \in \Acenter $.  Using these we facts, we get
	\begin{align*}
	 \widetilde{T_{\eta, \theta a}} & ( \omega  ) = \delta_{V_g ( \eta )} ( g (  \omega \tensora \theta a )  ) - \nabla_{V_g ( \eta )} V_g ( \theta a ) ( \omega )\\
	&= \delta_{V_g ( \eta )} ( g (  \omega \tensora \theta  ) a  ) - ( \nabla_{V_g ( \eta )} V_g ( \theta  ) a ) ( \omega ) - ( V_g ( \theta ) \delta_{V_g ( \eta )} ( a ) ) ( \omega ) ~ {\rm (} ~ {\rm by} ~{\rm Lemma } ~ \ref{9thnov20192}  {\rm )}\\
	&= \delta_{V_g ( \eta )} (  g ( \omega \tensora \theta )   ) a + g ( \omega \tensora \theta ) V_g ( \eta ) ( da ) - \nabla_{V_g ( \eta )} V_g ( \theta ) ( a \omega ) - ( V_g ( \theta ) g ( \eta \tensora da )   ) ( \omega )\\
	&= ( \delta_{V_g ( \eta )} ( g ( \omega \tensora \theta ) ) - \nabla_{V_g ( \eta )} V_g ( \theta ) ( \omega )  ) a + g ( \omega \tensora \theta ) g ( \eta \tensora da ) - g ( \eta \tensora da ) g ( \theta \tensora \omega )\\
	& \qquad\qquad {\rm (} ~ {\rm as} ~ a ~ \in ~ \Acenter ~ {\rm and} ~ V_g ( \theta )  \in  \mathcal{X} ( \A )  =  \mathcal{Z} ( \E^* ) ~ {\rm )} \\
	&=  \widetilde{T_{\eta, \theta }} ( \omega  ) a.
	\end{align*}
	This proves the lemma.
	\end{proof}
	Let us then define, for $ \omega \in \E $ the map
	$$ \widetilde{T} ( \omega ): \Ecenter \otimes_{\Acenter} \Ecenter \rightarrow \A, \qquad \widetilde{T} ( \omega ) ( \eta \otimes_{\Acenter} \theta ) = \widetilde{T_{\eta,\theta}} ( \omega ). $$
	It can be easily checked that $ \widetilde{T} ( \omega ) $ is well-defined
	and by Lemma \ref{13thapril2019lemma1}, it is right $ \Acenter $-linear. Moreover, if $ \omega \in \Ecenter$ one has $\widetilde{T_{\eta,\theta}} ( \omega ) \in \Acenter $ and thus $ \widetilde{T} ( \omega ) ( \eta \otimes_{\Acenter} \theta ) \in \Acenter $. Now, since we are taking $ \E = \Ecenter \otimes_{\Acenter} \A $, this allows us to make the identification 
	$$ \E \tensora \E = \Ecenter \otimes_{\Acenter} \Ecenter \otimes_{\Acenter} \A.$$
	Using this, it follows that  for $ \omega \in \Ecenter $, the map $ \widetilde{T} ( \omega ) $ extends to an element $ \widetilde{T_{{\rm ext}}} ( \omega ) $ of  $ \Hom_\A ( \E \tensora \E, \A ) $
	by the formula 
	$$  \widetilde{T_{{\rm ext}}} ( \omega \otimes_{\Acenter} \eta \otimes_{\Acenter} a ) = \widetilde{T} ( \omega \otimes_{\Acenter} \eta ) a;~ \omega, \eta \in \Ecenter,~ a \in \A.   $$

\noindent
We are ready for:	\\
{\bf Proof of Proposition \ref{15thapril2019prop}:} Uniqueness follows from Proposition \ref{15thapril2019prop2}. 
We need to prove the existence of the map $ \nabla $.  Let us recall (Proposition \ref{vg2nondegenerate}) that $ V_{g^{(2)}} $ is an isomorphism from $ \E \tensora \E $ to $ \Hom_\A ( \E \tensora \E, \A )$.    Hence, for $ \omega \in \Ecenter $, it makes sense to define 
	   $$ \nabla ( \omega ) = \vgtwo^{-1} \, (  \widetilde{T_{{\rm ext}}} ( \omega ) ). $$
	Next, for  $ \omega \in \Ecenter $ and $ a \in \A $, we define $ \nabla ( \omega a ) $ by the formula:
	$$ \vgtwo (  \nabla ( \omega a )  ) = \vgtwo ( \nabla ( \omega ) a + \omega \tensora da ). $$
	We claim that $ \nabla ( \omega a ) $ is well-defined. Indeed, if for   
	$ i = 1,2, \cdots n$ the elements $ \omega_i \in \Ecenter $ and $ a_i \in\A $ are such that $ \sum\nolimits_i \omega_i a_i = 0 $, then for all $ \eta, \theta \in \Ecenter, $
			\begin{align*}
			\sum\nolimits_i \vgtwo & ( \nabla ( \omega_i ) a_i + \omega_i \tensora d a_i  ) ( \eta \tensora \theta )\\
			&= \sum\nolimits_i [ \vgtwo ( \nabla ( \omega_i ) ) ( a_i \eta \tensora \theta ) + \vgtwo ( \omega_i \tensora d a_i ) ( \eta \tensora \theta )  ]\\
			&= \sum\nolimits_i [ \vgtwo ( \nabla ( \omega_i ) ) ( \eta \tensora \theta ) a_i + g^{(2)} ( ( \omega_i \tensora d a_i  ) \tensora ( \eta \tensora \theta )  ) ] ~ {\rm (} ~ {\rm as} ~ \eta,~ \theta \in \Ecenter ~ {\rm )} \\
			&= \sum\nolimits_i [ \delta_{V_g ( \eta )} ( g ( \theta \tensora \omega_i ) ) a_i - \nabla_{V_g ( \eta )} V_g ( \theta ) ( \omega_i ) a_i + g ( \omega_i \tensora g ( d a_i \tensora \eta ) \theta )  ]\\
			& \qquad\qquad \qquad {\rm (} ~ {\rm by} ~ {\rm the} ~ {\rm definitions} ~ {\rm of} ~ \widetilde{T} ( \omega ) ~ {\rm and} ~ \widetilde{T_{\eta,\theta}} ( \omega ) ~ {\rm )}\\
			&= \sum\nolimits_i [ \delta_{V_g ( \eta )} ( g ( \theta \tensora \omega_i ) a_i ) - g ( \theta \tensora \omega_i ) V_g ( \eta ) ( d a_i )  \\
			& \qquad\qquad \qquad - \nabla_{V_g ( \eta )} V_g ( \theta ) ( \omega_i a_i ) + g ( \omega_i \tensora \theta  ) g ( d a_i \tensora \eta )   ]\\
			& \qquad\qquad \qquad {\rm (} ~ {\rm as} ~ \delta ~ {\rm is} ~ {\rm a} ~ {\rm derivation} ~ {\rm and} ~ \theta \in \Ecenter ~ {\rm )}\\
			&= \delta_{V_g ( \eta )} ( g ( \theta \tensora \sum\nolimits_i \omega_i a_i ) ) - \sum\nolimits_i  g ( \theta \tensora \omega_i ) g ( \eta \tensora d a_i ) \\\ & \qquad\qquad \qquad 
			- \nabla_{V_g ( \eta )} V_g ( \theta ) ( \sum\nolimits_i \omega_i a_i  )
			+ \sum\nolimits_i g ( \omega_i \tensora \theta ) g ( d a_i \tensora \eta )\\
			&= \sum\nolimits_i g ( \omega_i \tensora \theta ) g ( d a_i \tensora \eta ) - \sum\nolimits_i g ( \theta \tensora \omega_i ) g ( \eta \tensora d a_i ) .
			\end{align*}
			In turn, as $ \eta, \theta \in \Ecenter $, we get $ g ( \omega_i \tensora \theta ) = g ( \theta \tensora \omega_i ) $ and $ g ( \eta \tensora da_i ) = g ( d a_i \tensora \eta )  $ and thus the above expression is equal to zero. Since $ \{ \eta \tensora \theta: \eta, \theta \in \Ecenter \} $ is right $\A$-total in $ \E \tensora \E $ (Proposition \ref{5thdec20191}), this proves that $ \nabla ( \omega a ) $ is well-defined.
			
			Since $ \E $ is centered, we have been able to extend the map $ \nabla $ to the whole of $ \E $.  Moreover, by the definition of $ \nabla $, for all $ \omega \in \Ecenter $ and $ a \in \A, $
	\begin{equation}	\label{15thapril20192}	 \nabla ( \omega a ) = \nabla ( \omega ) a + \omega \tensora da. \end{equation}
				It remains to show that $ \nabla $ is a connection. Since $ \E $ is centered, it suffices to prove the following equality for all $ \omega \in \Ecenter $ and for all $ a, b \in \A: $
		\begin{equation} \label{15thapril20193} \nabla ( ( \omega a  ) b ) = \nabla ( \omega a ) b + \omega a \tensora db. 	\end{equation}	
		But this follows by a simple computation using \eqref{15thapril20192}:
		\begin{align*}
		  \nabla ( \omega a b ) &= \nabla ( \omega ) a b + \omega \tensora d ( a b ) = \nabla ( \omega ) ab + \omega \tensora da b + \omega \tensora a db\\
		 &= \nabla ( \omega ) ab + \omega \tensora da b + \omega a \tensora  db = ( \nabla ( \omega ) a  + \omega \tensora da ) b + \omega a \tensora  db\\
		 &=  \nabla ( \omega a ) b + \omega a \tensora db.
		\end{align*}
		This finishes the proof of \eqref{15thapril20193} and hence the proposition.
	\qed
	
 \section{The Koszul formula in the covariant derivative form}	\label{section5}

Now we are ready to derive a Koszul formula for the Levi-Civita connection in the covariant derivative form (Theorem \ref{koszulformula1} and Theorem \ref{koszulformula2}). On the way to the derivation of these theorems, we will recover the classical formulas for the torsion and metric compatibility of a connection (Proposition \ref{torsionless26thjan} and Proposition \ref{metriccompatibility26thjan}). We start with a necessary and sufficient condition for a connection to be torsionless in the sense of Definition \ref{rLr}.

\subsection{The condition to be torsionless}

\bppsn \label{torsionless26thjan}
Suppose $ ( \E, d ) $ is a tame differential calculus and $\nabla$ is a connection on $\E$.    Then $ \nabla $ is  torsionless  if and only if for all $ X,  Y  \in  \mathcal{X} ( \A )$, 
		$$ \nabla_X Y - \nabla_Y X - [  X, Y ] = 0.  $$
	\eppsn
		\begin{proof} 
We will use Lemma 3.6 of \cite{article3} which
		states that if $\phi, \psi$ belong to $\E^*$ and $ W = \wedge \gamma \in \twoform $, for some $\gamma \in \E \tensora \E $, then, 
		\begin{equation} \label{8thdec20191}  ( \phi \tensora \psi  ) W = 2 ( \phi \tensora \psi ) ( 1 - \Psym ) ( \gamma ).  \end{equation}
		Now suppose $ \omega \in \E $.  With Sweedler  notation $ \nabla ( \omega ) = \omegazero \tensora \omegaone,  $ where $ \omegaone $ belongs to $\Ecenter$ (Remark \ref{sweedler}).
			Then we apply \eqref{8thdec20191}, using  $\Psym = \frac{1 + \sigma}{2}$ (Definition \ref{tame}), 
to compute 
					\begin{align*}
					 ( X \tensora Y ) \wedge \nabla ( \omega ) &= 2 ( X \tensora Y ) ( 1 - \Psym ) \nabla ( \omega )\\
				&=  ( X \tensora Y ) ( 1 - \sigma ) ( \omegazero \tensora \omegaone ) \\
				&= ( X \tensora Y ) ( \omegazero \tensora \omegaone - \omegaone \tensora \omegazero ).
				\end{align*}
				Now by using part 1. of Lemma \ref{lemma26thjan1}, we observe that
				$$ ( X \tensora Y ) ( \omegaone \tensora \omegazero ) = X ( \omegaone ) Y ( \omegazero ) = Y ( \omegazero ) X ( \omegaone ) = ( Y \tensora X ) \nabla ( \omega ). $$
				Therefore, for any connection $ \nabla $ on $ \E $ and $ X, Y $ in $ \mathcal{X} ( \A ) $, we have
		\begin{align} 
			( X \tensora Y ) \wedge \nabla ( \omega ) &= ( X \tensora Y ) \nabla ( \omega ) - ( Y \tensora X ) \nabla ( \omega ) \nonumber
			\\
		                                            &= \delta_Y ( X ( \omega ) ) - \delta_X ( Y ( \omega ) ) + ( \nabla_X Y - \nabla_Y X  ) ( \omega ) \label{29thjan2}
			\end{align}																					
	  by \eqref{30thjan3}. On the other hand, if $\nabla_0$ is the torsionless connection of Proposition \ref{5thdec20191}, then
				\begin{align}
				  - ( X \tensora Y ) d \omega \nonumber &= ( X \tensora Y ) \wedge \nabla_0 ( \omega )  \nonumber \\
				 &= ( X \tensora Y - Y \tensora X ) \nabla_0 ( \omega ) ~ 
				 \quad (\textup{by ~ the ~ first ~ equality in ~\eqref{29thjan2}}) \nonumber \\
				 &= - \delta_X ( Y ( \omega ) ) + \delta_Y ( X ( \omega ) ) + [ X, Y ] ( \omega ) \label{29thjan3}
				 \end{align}
				where in the last step, we have used Proposition \ref{prop26thjan2}.
				
		Next, suppose $ \nabla $ is a torsionless connection so that $ \wedge \nabla ( \omega ) = - d \omega $ and hence
				$$ ( X \tensora Y ) \wedge \nabla ( \omega ) = - ( X \tensora Y ) d \omega.  $$
				Comparing the second equality in \eqref{29thjan2} and \eqref{29thjan3}, we have
				$$ - \delta_X ( Y ( \omega ) ) + \delta_Y ( X ( \omega ) ) + ( \nabla_X Y - \nabla_Y X ) ( \omega ) = - \delta_X ( Y \omega ) + \delta_Y ( X ( \omega ) ) + [ X, Y ] ( \omega ) $$
				and therefore, for all $ \omega \in \E $, we deduce that
				$$ ( \nabla_X Y - \nabla_Y X ) ( \omega ) = [ X, Y ] ( \omega ). $$
				
\noindent				Conversely, suppose that for all $ X, Y  \in  \mathcal{X} ( \A )$, the equation $ \nabla_X Y - \nabla_Y X - [  X, Y ] = 0$ holds. Then by using the second equality in \eqref{29thjan2} and \eqref{29thjan3},  it is easy to check that for all $ X, Y \in \mathcal{X} ( \A ) $ and for all $ \omega \in \E $, 
				\begin{equation} \label{15thapril20194} ( X \tensora Y ) ( \wedge \circ \nabla ( \omega ) ) = - ( X \tensora Y ) ( d \omega ). \end{equation}
				If $ \nabla_0 $ is the torsionless connection of Proposition \ref{5thdec20191}, then $ \wedge \circ \nabla_0 ( \omega ) = - d \omega.  $ 		Hence, by virtue of \eqref{8thdec20191} and \eqref{15thapril20194},  we obtain for all $ X, Y $ in $ \mathcal{X} ( \A )$, 
				\begin{align*} 
				2 ( X \tensora Y ) ( ( 1 - \Psym  ) ( \nabla ( \omega ) - \nabla_0 ( \omega ) )   ) &= ( X \tensora Y ) \wedge \circ \nabla ( \omega ) - ( X \tensora Y ) \wedge \circ \nabla_0 ( \omega ) \\
				&=  - ( X \tensora Y ) d \omega + ( X \tensora Y ) d \omega = 0. 
				\end{align*}
			By an verbatim adaptation of the proof of Proposition \ref{15thapril2019prop2},	this allows us to conclude that 
			 $$	( 1 - \Psym  ) ( \nabla ( \omega ) - \nabla_0 ( \omega ) ) = 0. $$
			Applying $ \wedge $ to the equation and using the fact that $ {\rm Ran} ( \Psym ) = {\rm Ker} ( \wedge ) $ (Definition \ref{tame}), we conclude that
			$$ \wedge \circ \nabla ( \omega ) =\wedge  \nabla_0 ( \omega ) = - d \omega.$$
			Hence $\nabla $ is a torsionless connection.	This finishes the proof of the proposition.
				\end{proof}

				\subsection{The condition for metric-compatibility}  
		Next we come to proving a necessary and sufficient condition for a connection to be compatible with a metric $g$. For this we will need a couple of lemmas. Given a pseudo-Riemannian bilinear metric $g$ on $\E$,  we have a canonical  $\A$-$\A$-bilinear map $ \widetilde{g} : \E^* \tensora \E^* \rightarrow \A $ which we introduce in the next lemma. 
		\blmma \label{lemma26thjan2}
Let $ \widetilde{g}: \E^* \tensora \E^* \rightarrow \A $ be defined  by
 $$\widetilde{g} ( \phi \tensora \psi ) = g ( V^{-1}_g ( \phi ) \tensora V^{-1}_g ( \psi ) ).$$
It follows that for all  $ \phi, \psi \in \E^*, ~  $
 \begin{equation} \label{26thjan2} \phi ( V^{-1}_g ( \psi ) ) = \widetilde{g} ( \phi \tensora \psi ) \end{equation}
 and $\widetilde{g}$ is an $\A$-$\A$-bilinear map. 
\elmma
\begin{proof} The $\A$-bilinearity of $ \widetilde{g} $ follows from the bilinearity of $g$.     The equation \eqref{26thjan2} follows by a simple computation. Indeed,
  $$ \phi ( V^{-1}_g ( \psi ) ) = V_g ( V^{-1}_g ( \phi ) ) ( V^{-1}_g ( \psi ) ) = g ( V^{-1}_g ( \phi ) \tensora V^{-1}_g ( \psi ) ) = \widetilde{g} ( \phi \tensora \psi ) $$
 proving the lemma. \end{proof}

As a consequence, we have the following:
\blmma \label{lemma26thjan3}
	If $ \phi $ and $ \psi $ are elements of $ \E^* $ so that at least one of them is in $\mathcal{X} ( \A )$,  then 
	\begin{equation} \label{26thjan7} \widetilde{g} ( \phi \tensora \psi ) = \widetilde{g} ( \psi \tensora \phi ).  \end{equation}
	Moreover, if $X, Y \in \mathcal{X} ( \A )  $ and $ \omega \in \Ecenter $, then  
\begin{equation} \label{26thjan8} X ( \omega ) \in \mathcal{Z} ( \A ) \quad \textup{and} \quad \widetilde{g} ( X \tensora Y ) \in \mathcal{Z} ( \A ). \end{equation}
	\elmma
	\begin{proof} We use that $ \mathcal{X} ( \A ) = V_g ( \Ecenter ) $.  Then \eqref{26thjan7} follows by a combination of parts 1. and 2. of Lemma \ref{lemma2}. Similarly, both inclusions in \eqref{26thjan8} follows from part 3. of Lemma \ref{lemma2}. \end{proof}
	
	Before stating our main result, let us recall from condition \eqref{10thaugust20191} that a connection $\nabla$ is defined to be compatible with $g$ on the whole of $\E$ if $ \Pi_g ( \nabla ) ( e \tensora f ) = d ( g ( e \tensora f ) ) $ for all $e, f \in \E $, where the map $ \Pi_g ( \nabla ) $ is as defined in \eqref{22ndjuly20195}. 
				
  \bppsn \label{metriccompatibility26thjan}
	Suppose $ ( \E, d ) $ is a tame differential calculus.  A connection $ \nabla $ on $ \E $  is compatible on $ \E $ with the metric $g$ if and only if, for all $X,Y, Z$ in $ \mathcal{X} ( \A ) $, we have
		\begin{equation} \label{15thapril20195} \delta_Y ( \widetilde{g} ( Z \tensora X )  ) = \widetilde{g} ( \nabla_Y Z \tensora X  ) + \widetilde{g} (  \nabla_Y X \tensora Z ). \end{equation}
	\eppsn
	\begin{proof} Throughout the proof, we will use point 2. of Lemma \ref{lemma2}. Suppose $ \omega, \theta, \eta $ are unique elements in $ \Ecenter $ such that $ X = V_g ( \omega ), Y = V_g ( \theta )$ and $Z = V_g ( \eta ) $.  If $ \nabla $ is a connection on $ \E $ compatible with  $g$,  we get
	$$ Y \Big( ( g \tensora {\rm id} ) \big[ \sigma_{23} ( \nabla ( \omega ) \tensora \eta )  + ( \omega \tensora \nabla ( \eta ) \big] \Big) = Y ( dg ( \omega \tensora \eta ) ).$$
	 The left hand side of the above equation is equal to
	  \begin{align}
		Y \Big( ( g & \tensora {\rm id} ) \big[ \sigma_{23} ( \nabla ( \omega ) \tensora \eta )  + ( \omega \tensora \nabla ( \eta ) \big] \Big) \nonumber \\
		&= Y ( g ( \omegazero \tensora \eta ) \omegaone ) + Y ( g ( \omega \tensora \etazero  ) \etaone )  \nonumber \\
		&= g ( \omegazero \tensora \eta ) g ( \theta \tensora \omegaone ) + g ( \omega \tensora \etazero  ) g ( \theta \tensora \etaone ) \quad {\rm (}  ~ Y ~ {\rm being} ~ {\rm left} ~ \A-{\rm linear} ~  {\rm )} \nonumber \\
		&= ( V_g ( \eta ) \tensora V_g ( \theta )   ) \nabla ( \omega ) +  ( V_g ( \omega ) \tensora V_g ( \theta )   ) \nabla ( \eta )  \nonumber \\
		&= ( Z \tensora Y ) \nabla ( \omega ) + ( X \tensora Y ) \nabla ( \eta ) \nonumber \\
		&= \delta_Y ( Z ( \omega ) ) - ( \nabla_Y Z ) ( \omega ) + \delta_Y ( X ( \eta ) ) - (  \nabla_Y X ) ( \eta )  
		\qquad {\rm (} {\rm by} ~ \eqref{30thjan3}  {\rm )}  ~ \label{30thjan2} \, .
		\end{align}
	On the other hand, by using Lemma \ref{lemma26thjan3}, we obtain 
	\begin{equation} \label{10thdec20192} 
	Y ( dg ( \omega \tensora \eta ) )  = \delta_Y ( g ( V^{-1}_g ( X ) \tensora V^{-1}_g ( Z )  ) ) = \delta_Y ( \widetilde{g} ( X \tensora Z  )  ) = \delta_Y ( \widetilde{g} ( Z \tensora X  )  ). 
	\end{equation}
		 By combining \eqref{30thjan2} and \eqref{10thdec20192}, we have
		$$ \delta_Y ( Z ( \omega ) ) + \delta_Y ( X ( \eta ) ) - \delta_Y (  \widetilde{g} (  Z \tensora X  )  ) = (  \nabla_Y Z ) ( \omega ) + ( \nabla_Y X ) ( \eta ). $$
		However,
		\begin{align*}
		 \delta_Y ( Z ( \omega ) ) + \delta_Y ( X ( \eta ) ) & - \delta_Y (  \widetilde{g} (  Z \tensora X  )  ) \\ &= Y ( d g ( \eta \tensora \omega )  ) + Y (  d g ( \omega \tensora \eta )  ) - Y (  d g ( \eta \tensora \omega )  )\\
		&=  Y (  d g ( \eta \tensora \omega )  ) =  \delta_Y ( \widetilde{g} ( Z \tensora X  ) ).
		\end{align*}
		Therefore, using \eqref{26thjan2}
		$$ \delta_Y ( \widetilde{g} ( Z \tensora X  ) ) = ( \nabla_Y Z ) ( \omega ) + ( \nabla_Y X ) ( \eta ) = \widetilde{g} ( \nabla_Y Z \tensora X  ) + \widetilde{g} (  \nabla_Y X \tensora Z ) .
		$$	
		Conversely, suppose $ \nabla $ is a connection on $ \E $ such that equation \eqref{15thapril20195} is satisfied. We need to prove that $ \nabla $ is compatible with $g$ on the whole of $\E$.    Suppose $ \omega, \eta \in \Ecenter $ and $ Y \in \mathcal{X} ( \A ) $.  
		We claim that 
		$$ 
		Y \Big( ( g \tensora {\rm id} ) \big[\sigma_{23} (  \nabla ( \omega ) \tensora \eta ) + \omega \tensora \nabla ( \eta )   \big]   \Big) = Y ( dg ( \omega \tensora \eta ) ).
		$$
		Suppose $ X = V_g ( \omega ) $ and $ Z = V_g ( \eta ). $
		Then by applying \eqref{26thjan2} again, we get 
		\begin{align*}
		Y \Big( ( g \tensora & {\rm id} ) \big[\sigma_{23} (  \nabla ( \omega ) \tensora \eta ) + \omega \tensora \nabla ( \eta )   \big]   \Big)\\
		&= \delta_Y ( Z ( \omega )  ) - (  \nabla_Y Z ) ( \omega ) + \delta_Y ( X ( \eta ) ) - (  \nabla_Y X ) ( \eta ) ~ {\rm (} ~ {\rm by} ~ \eqref{30thjan2} ~ {\rm )} \\
		&= \delta_Y ( Z ( V^{-1}_g ( X )  ) ) - ( \nabla_Y Z ) (  V^{-1}_g ( X )  ) + \delta_Y ( X ( V^{-1}_g ( Z )  )  ) - ( \nabla_Y X  ) ( V^{-1}_g ( Z ) )\\
		&= \delta_Y ( \widetilde{g} ( Z \tensora X  )  ) - \widetilde{g} ( \nabla_Y Z \tensora X ) + \delta_Y (  \widetilde{g} (  X \tensora Z  )  ) - \widetilde{g} (  \nabla_Y X \tensora Z  )\\
		&= \widetilde{g} ( \nabla_Y Z \tensora X ) + \widetilde{g} (  \nabla_Y X \tensora Z  ) ~ {\rm (} ~ {\rm by} ~ \eqref{15thapril20195} ~ {\rm )}\\
	  &= \delta_Y (  \widetilde{g} ( Z \tensora X  )  ) \qquad {\rm (} ~ {\rm by} ~ \eqref{15thapril20195} ~ {\rm )}\\
		&= Y ( dg ( \omega \tensora \eta ) ) \qquad {\rm (} ~ {\rm by} ~ \eqref{10thdec20192} ~ {\rm )}.
		\end{align*}
	This proves the claim and hence,
		$$ ( \Pi_g ( \nabla ) - dg ) ( \omega \tensora \eta ) = 0 ~ {\rm for} ~ {\rm all} ~ \omega, \eta \in \Ecenter.$$
		Here, the map $ \Pi_g $ is as in \eqref{22ndjuly20195}. However, by  Proposition 4.7 of \cite{article1}, the map $ \Pi_g ( \nabla ) - dg $ is right $\A$-linear. Since the set $ \{ \omega \tensora \eta: \omega, \eta \in \Ecenter \} $ is right $\A$-total in $ \E \tensora \E $, one has $ \Pi_g ( \nabla ) - dg = 0 $ and so $ \nabla $ is compatible with $g$ on the whole of $\E$.    
		\end{proof}
		
		\subsection{The Koszul formula}
		Now we are in a position to spell out the Koszul formula of the Levi-Civita connection in the covariant derivative form. We give a sketch of the proofs since these involve techniques already employed in this article. We will need the following simple corollary to Lemma \ref{lemma26thjan3}.
		
\blmma \label{lemma27thjan1}
Suppose $X, Y, Z \in \mathcal{X} ( \A ) $ and $ a \in \A $.  Then the following formulas hold:
\begin{align*}
\widetilde{g} ( Z \tensora [ X, Y ] ) & = \widetilde{g} ( [ X, Y ] \tensora Z ), \\  \delta_Z (  \widetilde{g} ( X \tensora Y ) ) a & = \delta_{Za} (  \widetilde{g} ( X \tensora Y ) ).
\end{align*}
\elmma
\begin{proof} The first equality follows from \eqref{26thjan7}. For the second, we claim $ \delta_{Z} ( \widetilde{g} ( X \tensora Y  ) ) \in \Acenter $.  Indeed, by \eqref{26thjan8}, $ \widetilde{g} ( X \tensora Y ) \in \Acenter $ and so by part 4. of Lemma \ref{lemma2},  $ d ( \widetilde{g} ( X \tensora Y  )   ) \in \Ecenter $.  Since $ Z \in \mathcal{X} ( \A ) $, there exists a unique $ \eta \in \Ecenter $ such that $ V_g ( \eta ) = Z $ and thus by Lemma \ref{lemma2}, we can conclude that
$$ \delta_Z ( \widetilde{g} (  X \tensora Y  ) ) = g ( \eta \tensora d ( \widetilde{g} (  X \tensora Y )  ) ) ~ \in ~ \Acenter. $$
This proves the claim. Hence, since $ Z $ belongs to $ \mathcal{X} ( \A ) = \mathcal{Z} ( \E^* ) $, we have
\begin{align*}
\delta_Z (  \widetilde{g} ( X \tensora Y ) ) a & = a ( \delta_Z (  \widetilde{g} ( X \tensora Y ) ) ) = a Z ( d ( \widetilde{g} ( X \tensora Y ) ) )\\
 &= ( Z a ) ( d ( \widetilde{g} ( X \tensora Y ) ) )  = \delta_{Za} (  \widetilde{g} ( X \tensora Y ) ).
\end{align*}
This finishes the proof of the lemma.  \end{proof}
		
		\bppsn \label{koszulformula1}
	Suppose $ ( \E, d ) $ is a tame differential calculus and $g$ is a pseudo-Riemannian bilinear metric on $\E$.    If $ \nabla $ is the Levi-Civita connection for the pair $ ( \E, g ) $ and $ X, Y, Z $ are elements of $\mathcal{X} ( \A )$,   then we have the following Koszul formula for $\nabla:$ 
	\begin{align*} 2 \widetilde{g} ( \nabla_X Y \tensora Z ) 
	 & = \delta_X ( \widetilde{g} ( Y \tensora Z  ) )  + \delta_Y ( \widetilde{g} ( X \tensora Z  ) )
	- \delta_{Z} ( \widetilde{g} ( X \tensora Y ) )    \\  & \quad - \widetilde{g} ( Y \tensora [ X, Z  ] )  - \widetilde{g} (  [ Y , X ] \tensora Z ) + \widetilde{g} ( X \tensora [ Z, Y ] ).
	\end{align*}
	\eppsn
	\begin{proof} We use Lemma \ref{lemma26thjan3}. We start with Proposition \ref{metriccompatibility26thjan} and permute $X,Y,Z$ cyclically. Then we employ  Proposition \ref{torsionless26thjan} and follow the proof  of Theorem 5.5 of \cite{article3} to derive the desired equation.
	\end{proof}
	
	Interestingly, the Koszul formula in Proposition \ref{koszulformula1} needs only $ X, Y  $ in $  \mathcal{X} ( \A ) $  while $ Z $ can be in $ \E^* $. 
	
	\bppsn \label{koszulformula2}
Under the hypotheses of Theorem \ref{koszulformula1}, if $ X, Y $ are elements in $ \mathcal{X} ( \A ) $ and $ Z \in \E^* $, then we have the following equality:
	\begin{align*}
	2 \widetilde{g} ( \nabla_X Y \tensora Z  ) & = 
	 \delta_X ( \widetilde{g} ( Y \tensora Z  ) )  + \delta_Y ( \widetilde{g} ( X \tensora Z  ) )  - \delta_{Z} ( \widetilde{g} ( X \tensora Y ) ) \\ 
	& \quad - \widetilde{g} ( Y \tensora [ X, Z  ] )  - \widetilde{g} (  [ Y , X ] \tensora Z ) + \widetilde{g} ( X \tensora [ Z, Y ] ).
	 \end{align*}
	Here, $ [ X, Z ] $ and $ [ Y, Z ] $ are defined as in \eqref{4thfeb20192} which is well-defined by Proposition \ref{4thfeb2019prop}.
	\eppsn
	\begin{proof} As $ \mathcal{X} ( \A ) $ is right $\A$-total in $\E^*$,   it suffices to prove the formula for  $ 2 \widetilde{g} ( \nabla_X Y \tensora Z a ) $ where $ Z \in \mathcal{X} ( \A ) $ and $ a \in \A $.  We use  the two equations of Lemma \ref{lemma27thjan1} and Proposition \ref{koszulformula1} to see that
	\begin{align*}
	 2 \widetilde{g} ( \nabla_X Y \tensora Z a ) & = 2 \widetilde{g} ( \nabla_X Y \tensora Z  ) a\\
	 &= \delta_X ( \widetilde{g} ( Y \tensora Z a ) )  + \delta_Y ( \widetilde{g} ( X \tensora Z a ) )  
	 - \delta_{Za} ( \widetilde{g} ( X \tensora Y ) ) \\  
	 & \quad - \widetilde{g} ( Y \tensora [ X, Z a ] )  - \widetilde{g} (  [ Y , X ] \tensora Za )
	- \widetilde{g} ( X \tensora [ Y, Z a ] ).
	\end{align*}
	This finishes the proof of the proposition.
	\end{proof}
	
	\vspace{4mm}
	
	\noindent
{\bf Acknowledgment:} JB and DG were funded by a ``Research in Pairs" INDAM grant at ICTP, Trieste, and by a ``Dipartimento di Eccellenza" fund of the Department of Mathematics and Geosciences, University of Trieste. DG was partially supported by JC Bose Fellowship and Grant from the Department of Science and Technology, Govt. of India. GL acknowledges partial support from INFN, Iniziativa Specifica GAST and from INDAM GNSAGA. Finally, JB and DG want to thank GL for the very kind hospitality during their stay in Trieste.

\vspace{1in}

\noindent
Jyotishman Bhowmick, Debashish Goswami: \\
Indian Statistical Institute, 
203, B. T. Road, Kolkata 700108, India\\
 jyotishmanb$@$gmail.com, goswamid$@$isical.ac.in\\

\noindent 
Giovanni Landi: \\ Matematica, Universit\`a di Trieste, Via A. Valerio, 12/1, 34127 Trieste, Italy \\
Institute for Geometry and Physics (IGAP) Trieste, Italy and INFN, Trieste, Italy \\
landi@units.it

\end{document}